# Avaliação computacional de esquemas convectivos em problemas de dinâmica dos fluidos

# Computational evaluation of convection schemes in fluid dynamics problems

Valdemir Garcia Ferreira[1]; Giseli Aparecida Braz de Lima[2]; Laís Corrêa[2]; Miguel Antonio Caro Candezano[2], Eliandro Rodrigues Cirilo[3]; Paulo Laerte Natti[3] e Neyva Maria Lopes Romeiro[3]

## Resumo

Este artigo apresenta uma avaliação computacional dos esquemas populares convectivos *upwind* de alta resolução WACEB, CUBISTA e ADBQUICKEST para resolver numericamente problemas não lineares em dinâmica dos fluidos. Os esquemas são analisados e implementados em variáveis normalizadas de Leonard e no contexto do método das diferenças finitas. O desempenho dos esquemas é avaliado em problemas de Riemann para as equações 1D de Burgers, de Euler e de águas rasas. Com base nos resultados numéricos, os esquemas são classificados segundo seus desempenhos em resolver essas equações não lineares. O esquema selecionado é, então, aplicado na simulação numérica de escoamentos incompressíveis 3D com superfícies livres móveis.

**Palavras-chave**: Esquemas de alta resolução. Simulação numérica. Termos convectivos. Leis de conservação. Escoamentos com superfícies livres.

## Abstract

This article provides a computational evaluation of the popular high resolution upwind WACEB, CUBISTA and ADBQUICKEST schemes for solving non-linear fluid dynamics problems. By using the finite difference methodology, the schemes are analyzed and implemented in the context of normalized variables of Leonard. In order to access the performance of the schemes, Riemann problems for 1D Burgers, Euler and shallow water equations are considered. From the numerical results, the schemes are ranked according to their performance in solving these non-linear equations. The best scheme is then applied in the numerical simulation of tridimensional incompressible moving free surface flows.

**Key words:** High resolution schemes. Numerical simulation. Convection terms. Conservation laws. Flow with free surfaces.

[1] Docente do Departamento de Matemática Aplicada e Estatística do Instituto de Ciências Matemática e de Computação da Universidade de São Paulo em São Carlos (ICMC-USP); pvgf@icmc.usp.br.

[2] Alunos do curso de doutorado em Ciências da Computação do ICMSC-USP/São Carlos; {giabl,lacorrea,macc}@icmc.usp.br.

[3] Alunos do curso de doutorado em Ciências da Computação do ICMSC-USP/São Carlos; {giabl,lacorrea,macc}@icmc.usp.br.

[4] Alunos do curso de doutorado em Ciências da Computação do ICMSC-USP/São Carlos; {giabl,lacorrea,macc}@icmc.usp.br.

[5] Docentes do Departamento de Matemática da Universidade Estadual de Londrina; {ercirilo,plnatti,nromeiro}@uel.br.

[6] Docentes do Departamento de Matemática da Universidade Estadual de Londrina; {ercirilo,plnatti,nromeiro}@uel.br.

[7] Docentes do Departamento de Matemática da Universidade Estadual de Londrina; {ercirilo,plnatti,nromeiro}@uel.br.

## Introdução

Aproximação de termos convectivos (em geral não lineares) em leis de conservação por esquemas *upwind* de alta resolução é uma área em constante atividade e tem atraído muitos pesquisadores em dinâmica dos fluidos computacional. Na literatura especializada, há uma variedade de esquemas para a representação desses termos; entretanto nenhum deles tem mostrado ser completamente robusto (WATERSON; DECONINCK, 2007). Como exemplo, pode-se citar os esquemas populares *Weighted-Average Coefficient Ensuring Boundedness* (WACEB) (SONG et al., 2000), *Convergent and Universally Bounded Interpolation Scheme for Treatment of Advection* (CUBISTA) (ALVES; OLIVEIRA; PINHO, 2003) e, mais recentemente, o novo esquema adaptativo *Quadratic Upstream for Convective Kinematics with Estimated Streaming Terms* (QUICKEST) (FERREIRA et al., 2007), também conhecido como ADBQUICKEST. Estes esquemas foram selecionados nesse trabalho por representarem bem o estado da arte atual em esquemas *upwind* de alta resolução para termos convectivos.

Dentro desse cenário, o objetivo do presente estudo consiste investigar o desempenho dos esquemas WACEB, CUBISTA e ADBQUICKEST, no contexto de variáveis normalizadas e utilizando-se a metodologia de diferenças finitas. A seguir são apresentados os três esquemas *upwind* considerados no trabalho.

- WACEB (SONG et al., 2000):

$$\hat{\phi}_f = \begin{cases} 2\hat{\phi}_U, & 0 \leq \hat{\phi}_U < \dfrac{3}{10}, \\ \dfrac{3}{4}\hat{\phi}_U + \dfrac{3}{8}, & \dfrac{3}{10} \leq \hat{\phi}_U \leq \dfrac{5}{6}, \\ 1, & \dfrac{5}{6} < \hat{\phi}_U \leq 1, \\ \hat{\phi}_U, & \text{caso contrário.} \end{cases} \quad (1)$$

- CUBISTA (ALVES; OLIVEIRA; PINHO, 2003):

$$\hat{\phi}_f = \begin{cases} \dfrac{7}{4}\hat{\phi}_U, & 0 < \hat{\phi}_U < \dfrac{3}{8}, \\ \dfrac{3}{4}\hat{\phi}_U + \dfrac{3}{8}, & \dfrac{3}{8} \leq \hat{\phi}_U \leq \dfrac{3}{4}, \\ \dfrac{1}{4}\hat{\phi}_U + \dfrac{3}{4}, & \dfrac{3}{4} < \hat{\phi}_U < 1, \\ \hat{\phi}_U, & \text{caso contrário.} \end{cases} \quad (2)$$

- ADBQUICKEST (FERREIRA et al., 2007):

$$\hat{\phi}_f = \begin{cases} (2-\theta)\hat{\phi}_U, & 0 < \hat{\phi}_U < a, \\ \hat{\phi}_U + \dfrac{1}{2}\tau(1-\hat{\phi}_U) - \dfrac{1}{6}\mu(1-2\hat{\phi}_U), & a \leq \hat{\phi}_U \leq b, \\ 1-\theta+\theta\hat{\phi}_U, & b < \hat{\phi}_U < 1, \\ \hat{\phi}_U, & \text{caso contrário.} \end{cases} \quad (3)$$

Em (3) tem-se que $\tau = (1-|\theta|)$, $a = \dfrac{2-3|\theta|+\theta^2}{7-6\theta-3|\theta|+2\theta^2}$,

$b = \dfrac{-4+6\theta-3|\theta|+\theta^2}{-5+6\theta-3|\theta|+2\theta^2}$ e $\mu = (1-\theta^2)$, onde $\theta = a\,\delta_t/\delta_x$

é o número de Courant, com $a$ uma velocidade de convecção, $\delta_t$ a marcha no tempo e $\delta_x$ o espaçamento espacial.

Nos esquemas acima, $\hat{\phi}_U$ é a variável normalizada de Leonard (LEONARD, 1988) definida por

$$\hat{\phi}_U = \dfrac{\phi_U - \phi_D}{\phi_D - \phi_R}, \quad (4)$$

onde $D$, $U$ e $R$ correspondem, com respeito ao fluxo de massa numa face $f$ da célula computacional, às posições à jusante, à montante e mais à montante, respectivamente (FERREIRA et al., 2007).

Este trabalho é apresentado na seguinte seqüência. Na seção 2 o desempenho dos esquemas WACEB, CUBISTA e ADBQUICKEST é avaliado em três problemas de Riemann para as equações não lineares 1D de Burgers, de Euler e de águas rasas. Na seção 3, o esquema ADBQUICKEST, que apresentou o melhor desempenho na seção precedente, é utilizado na solução numérica das equações 3D de Navier-Stokes, caso incompressível com superfícies livres móveis, para simular os fenômenos do salto hidráulico circular e de dobras em jatos livres a baixo número de Reynolds. Enfim, na seção 4, as conclusões deste trabalho são apresentadas.

## Experimentos Numéricos

Nesta seção, os esquemas WACEB, CUBISTA e ADBQUICKEST são usados na solução de problemas de Riemann para as equações de Burgers, Euler e águas rasas.

### Equação de Burgers

A equação não linear de Burgers

$$u_t + \left(\frac{u^2}{2}\right)_x = vu_{xx}, \quad 0 \le x \le 1, \quad t > 0, \quad (5)$$

onde $v$ é o coeficiente de viscosidade, modela o transporte da propriedade $u$, com velocidade de convecção u, dada a condição inicial $u(x,0) = u_0$.

Na seqüência são apresentadas as simulações de dois problemas teste: teste-1 e teste-2. No teste-1, a Eq.(5) é resolvida com a condição inicial $u(x,0)=sen(x)$, com $x \in [-1,1]$, para o coeficiente de viscosidade $v = 0 m^2 s^{-1}$, enquanto no teste-2, a Eq.(5) é resolvida com a condição inicial $u(x,0)= -sen(\pi x)$, com $x \in [-1,1]$, para dois coeficientes de viscosidade: $v = 0.1 m^2 s^{-1}$ e $v = 0.0001 m^2 s^{-1}$.

No teste-1, a solução exata é dada por (PLATZMAN, 1964):

$$u(x,t) = -2\sum_{n}^{\infty} \frac{J_n(-nt)}{nt} sen(nx), \quad (6)$$

em que $J_n$ é a função de Bessel de ordem $n$ de primeira espécie. Neste trabalho, considera-se a solução (6) truncada em termos $n=500$ (semi-analítica). Na simulação que segue, utilizaram-se de uma malha uniforme de 400 células computacionais, com tempo final de simulação equivalente a $t=1$ (instante que ocorre o choque) e número de Courant $\theta=0.5$ Os resultados numéricos usando os esquemas WACEB, CUBISTA e ADBQUICKEST estão apresentados na Figura 1. Nesta figura, a coluna à esquerda corresponde às comparações entre as soluções computadas e a solução semi-analítica, enquanto a coluna à direita corresponde aos erros absolutos.

**Figura 1 -** Soluções semi-analítica e numéricas (à esquerda) e erros absolutos (à direita), em escala log-log, gerados pelos esquemas WACEB (a)-(b), CUBISTA (c)-(d) e ADBQUICKEST (e)-(f).

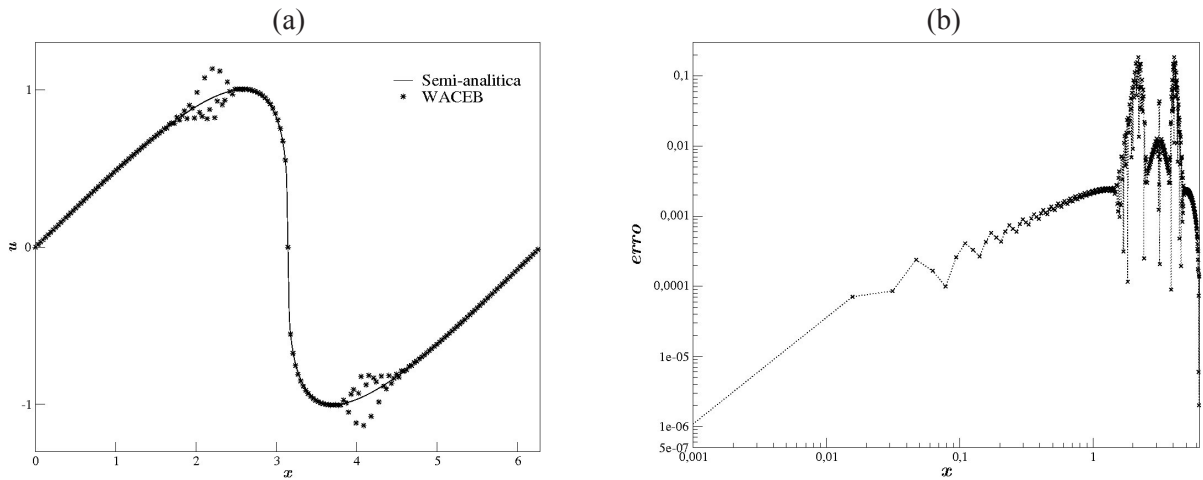

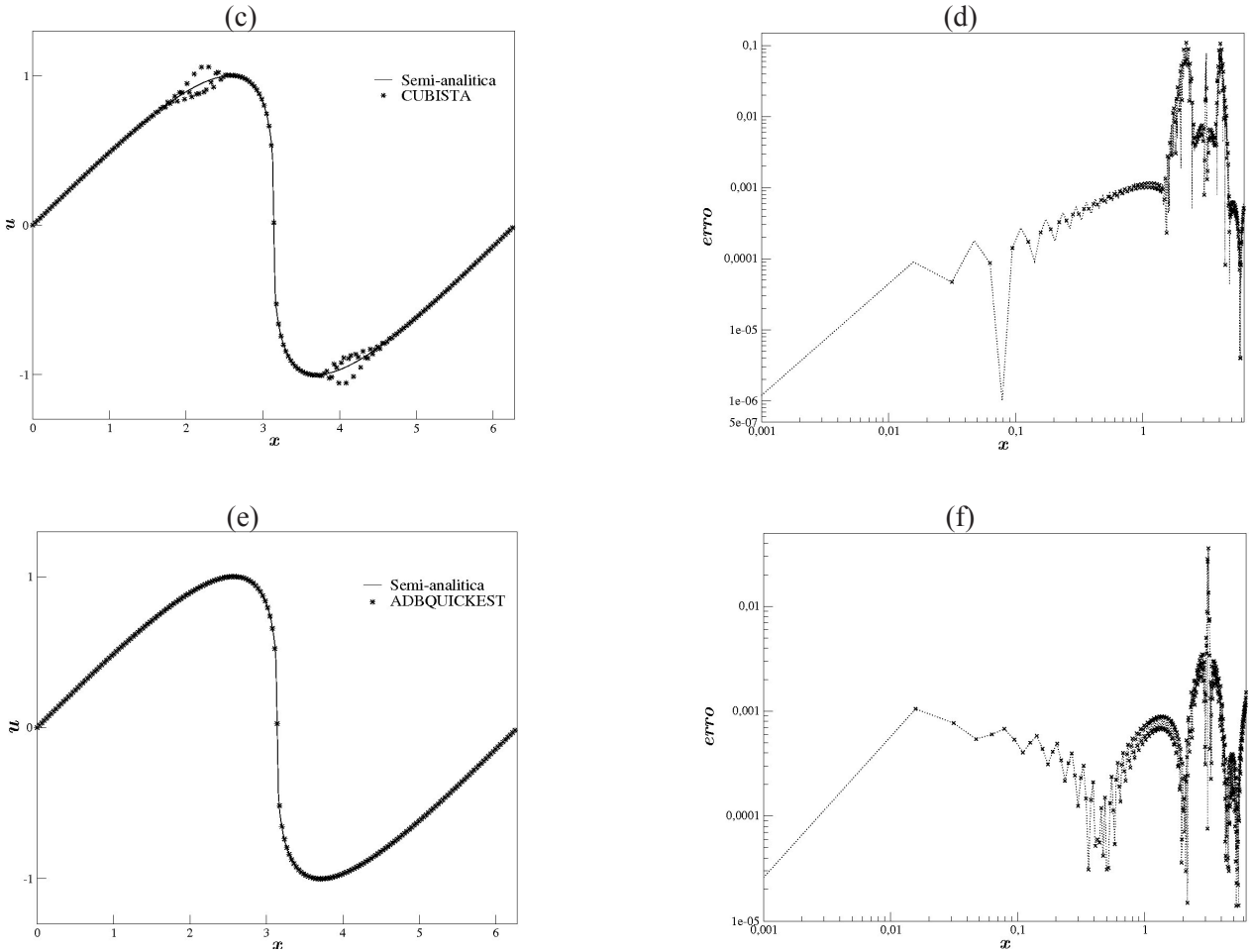

Na Figura 1 é possível observar que os melhores resultados foram obtidos com o esquema ADBQUICKEST. Note que os resultados dos outros esquemas apresentam oscilações evidentes. Constata-se ainda que os erros obtidos com os esquemas WACEB e CUBISTA são maiores que os erros cometidos no cálculo com o ADBQUICKEST. Com o objetivo de estimar a ordem de convergência dos métodos numéricos aqui apresentados, quando aplicados a estes problemas não lineares, considera-se a equação (5) sem viscosidade e com a condição inicial $u(x,0) = sen(x)$, $x \in [0,2\pi]$. Para tanto obtém-se o erro relativo $E_{\delta_x}$ numa malha de espaçamento $\delta_x$, com a estimativa para a ordem de convergência dada por $\tilde{p} = \log\left(\left\|E_{\delta_x}\right\|/\left\|E_{\delta_x/2}\right\|\right)/\log 2$, onde as normas

$L_1$, $L_2$ e $L_\infty$ são consideradas. Nas simulações, foram usadas malhas com 10, 20, 40, 80, 160 células computacionais e $\Theta = 0.2$. As estimativas para ordem $\tilde{p}$ são apresentadas na Tabela 1, onde pode-se observar que o método ADBQUICKEST é mais preciso.

**Tabela 1 -** Erros nas normas $L_1$, $L_2$ e $L_\infty$ e estimativas para ordem de convergência $\tilde{p}$.

| Esquemas | N | $L_1\,E_{\delta_x}$ | $\tilde{p}$ | $L_2\,E_{\delta_x}$ | $\tilde{p}$ | $L_\infty E_{\delta_x}$ | $\tilde{p}$ |
|---|---|---|---|---|---|---|---|
| WACEB | 10 | 0.157e-1 | — | 0.180e-1 | — | 0.214e-1 | — |
| | 20 | 0.708e-2 | 1.15 | 0.835e-2 | 1.11 | 0.130e-1 | 0.72 |
| | 40 | 0.413e-2 | 0.78 | 0.452e-2 | 0.89 | 0.695e-2 | 0.90 |
| | 80 | 0.224e-2 | 0.88 | 0.242e-2 | 0.90 | 0.479e-2 | 0.54 |
| | 160 | 0.117e-2 | 0.94 | 0.127e-2 | 0.93 | 0.269e-2 | 0.83 |
| CUBISTA | 10 | 0.199e-1 | — | 0.206e-1 | — | 0.233e-1 | — |
| | 20 | 0.796e-2 | 1.33 | 0.854e-2 | 1.27 | 0.100e-1 | 1.22 |
| | 40 | 0.361e-2 | 1.14 | 0.412e-2 | 1.05 | 0.577e-2 | 0.80 |
| | 80 | 0.168e-2 | 1.11 | 0.187e-2 | 1.14 | 0.350e-2 | 0.72 |
| | 160 | 0.800e-3 | 1.07 | 0.870e-3 | 1.11 | 0.172e-2 | 1.03 |
| ADBQUICKEST | 10 | 0.178e-1 | — | 0.182e-1 | — | 0.185e-1 | — |
| | 20 | 0.553e-2 | 1.69 | 0.614e-2 | 1.56 | 0.722e-2 | 1.36 |
| | 40 | 0.187e-2 | 1.57 | 0.232e-2 | 1.41 | 0.334e-2 | 1.11 |
| | 80 | 0.888e-3 | 1.07 | 0.101e-2 | 1.20 | 0.152e-2 | 1.14 |
| | 160 | 0.438e-3 | 1.02 | 0.473e-3 | 1.09 | 0.769e-3 | 0.98 |

O teste-2 é um modelo interessante que combina efeitos de advecção (não lineares) e efeitos viscosos (lineares), sendo, portanto, essencial no estudo da dinâmica dos fluidos. Ele serve também como um modelo simplificado para o entendimento da formação do choque e da turbulência. Para a simulação deste teste com o esquema ADBQUICKEST foram considerados $\Theta = 0.5$ e cinco tempos diferentes, a saber, $t = 0.1,\ 0.3,\ 0.5,\ 0.7,\ 0.9$. A Figura 2 mostra as soluções numéricas para os casos em que $v = 0.1\,m^2s^{-1}$ (Fig.2-(a)) e $v = 0.0001\,m^2s^{-1}$ (Fig.2-(b)). Em particular, observa-se a formação do choque em e a capacidade do esquema ADBQUICKEST em reproduzir essa descontinuidade sem qualquer vestígio de oscilações numéricas. Analogamente à situação observada no problema teste-1, os outros esquemas apresentaram oscilações nas vizinhanças dessa descontinuidade.

**Figura 2 -** Soluções numéricas em vários tempos para a equação 1D de Burgers com o esquema ADBQUICKEST:  (a) $v=0.1m^2s^{-1}$ e (b) $v=0.0001m^2s^{-1}$

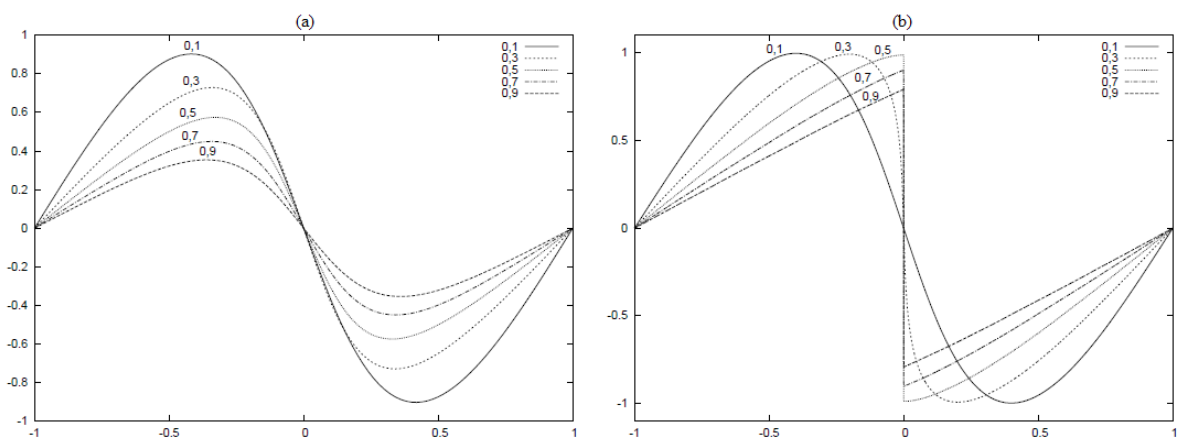

## Equações de Euler

Nesta subseção, simula-se um problema de Riemann, conhecido como tubo de choque (TORO, 2001), para as equações de Euler da dinâmica dos gases

$$\mathbf{U}_t + \mathbf{F(U)}_x = 0 \qquad (7)$$

em que $\mathbf{U} = \left[\rho, \rho u, E\right]^T$  e  $\mathbf{F} = \left[\rho u, \rho u^2 + p, u(E + p)\right]^T$. Aqui, as variáveis conservadas $\rho, u, \rho u, E, p$ são, respectivamente, a massa específica, a velocidade, a quantidade de movimento, a energia total e a pressão do sistema de partículas. Para a resolução do sistema de equações diferenciais parciais (7) é necessário uma equação de estado, assim utiliza-se a equação de gás ideal $p = (\gamma - 1)\left(E - 0.5\rho u^2\right)$, com $\gamma = 1.4$. A condição inicial proposta por (TITAREV; TORO, 2005)

$$(\rho, \upsilon, p)^T = \begin{cases} (1.515695,\ 0.523346,\ 1.80500)^T, & x < -4.5 \quad (8) \\ (1 + 0.1sen(20\pi x),\ 0.0,\ 1)^T, & x \geq -4.5 \end{cases}$$

foi considerada no domínio de cálculo $x \in [-5,5]$. Para a simulação desse problema, adotaram-se $\Theta = 0.6$ e uma malha de 12500 células computacionais. A solução de referência foi gerada por um esquema de primeira ordem, descrito em (SPALDING, 1972), para uma malha de 25000 células. A Figura 3 mostra as soluções de referência e as numéricas para a massa específica. Pode-se constatar novamente, por meio das ampliações locais da Figura 3, que o esquema ADBQUICKEST (caso (c)) produziu um resultado melhor que aqueles produzidos pelos esquemas WACEB (caso (a)) e CUBISTA (caso (b)).

## Equações de águas rasas

Equações de águas rasas são as equações (7) com $\mathbf{U} = \left[h, hu\right]^T$ e $\mathbf{F} = \left[hu, hu^2 + \frac{1}{2}gh^2\right]^T$, sendo $h$ a altura do fluido (no caso água), $u$ a velocidade horizontal da água num canal e  a aceleração devido à gravidade. O fenômeno modelado por este sistema hiperbólico é conhecido como ruptura de barragem (TORO, 2001). As simulações usando os três esquemas WACEB, CUBISTA e ADBQUICKEST numa malha de 400 células computacionais e $\Theta = 0.6$ estão apresentadas na Figura 4. A Figura 4(a) corresponde às comparações das soluções numéricas e de referência em todo domínio, enquanto as Figuras 4(b) e 4(c) correspondem às ampliações das regiões destacadas por retângulos na Figura 4(a). A solução de referência corresponde à calculada pelo método de Godunov numa malha de 10000 células computacionais. Mais uma vez, verifica-se que o esquema ADBQUICKEST forneceu os melhores resultados.

**Figura 3 -** Comparação das soluções numéricas e de referência para a massa específica nas equações de Euler: (a) WACEB; (b) CUBISTA e (c) ADBQUICKEST.

**Figura 4 -** (a) Comparações das soluções numéricas e de referência para *h* nas equações de águas rasas; (b) e (c) ampliações nos intervalos [-3,-2] e [1.6; 2.8]

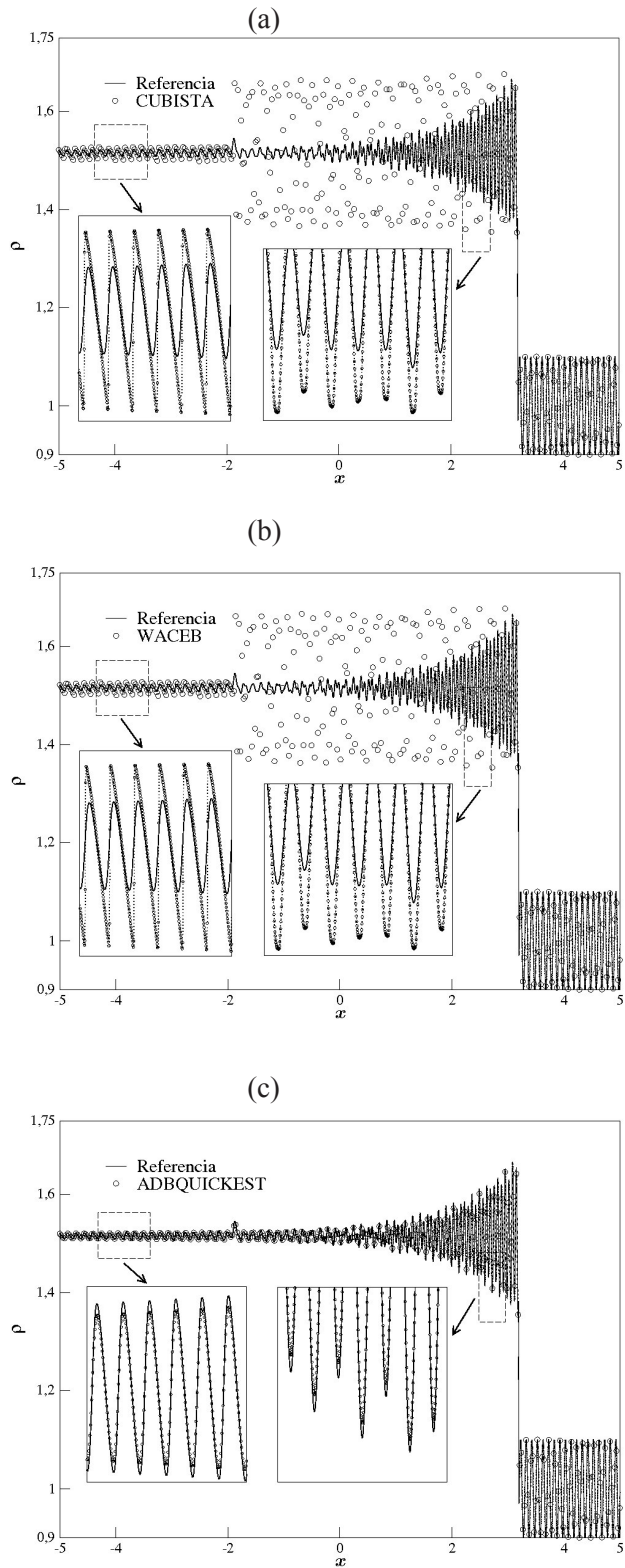

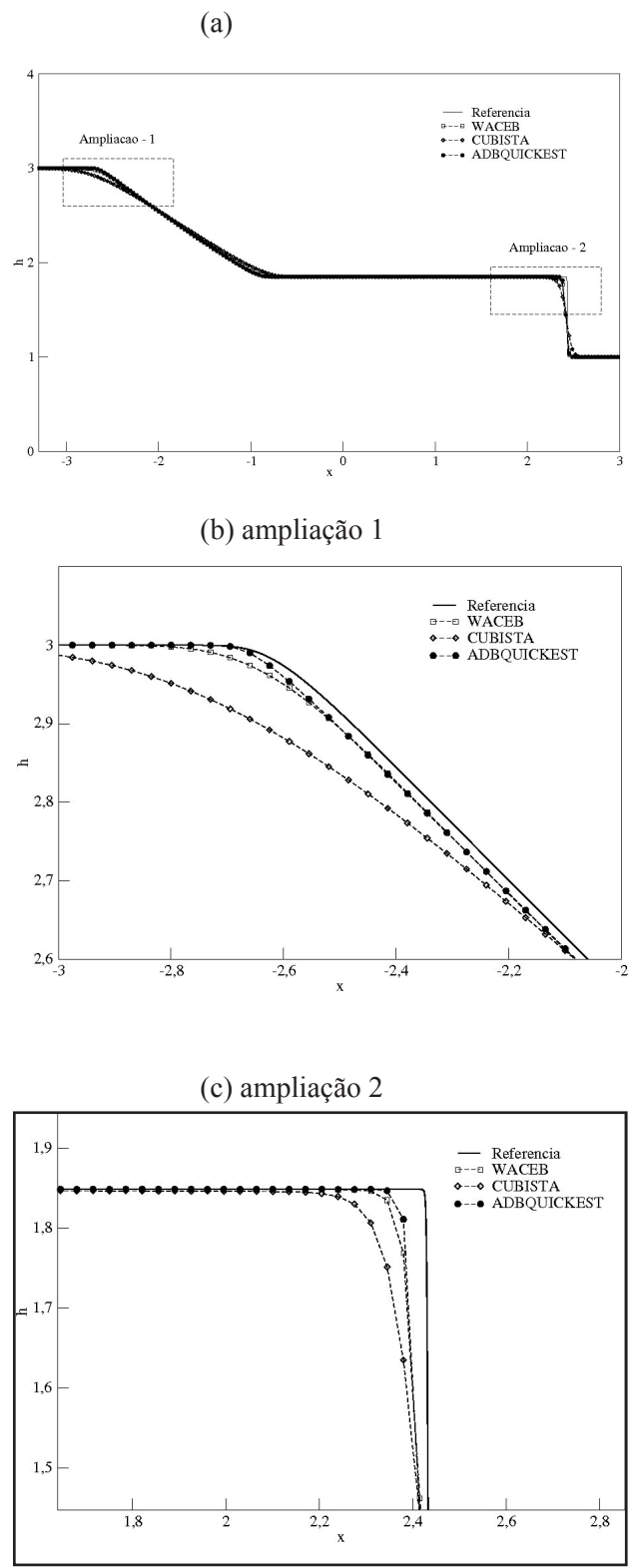

## Esquema ADBQUICKEST: escoamentos incompressíveis 3D

Verificado o bom desempenho do esquema ADBQUICKEST nos testes precedentes, apresentam-se agora simulações, com este esquema, de escoamentos incompressíveis tridimensionais com superfícies livres móveis.

As simulação realizadas são do salto hidráulico circular (ELLEGAARD et al., 1998) e de dobras em jatos livres a baixo número de Reynolds (TOMÉ et al., 2001). Tais escoamentos são modelados pelas equações adimensionais completas de Navier-Stokes, ou seja,

$$\frac{\partial v_i}{\partial t} + \frac{\partial (v_i v_j)}{\partial x_j} = -\frac{\partial p}{\partial x_i} + \frac{1}{R_e}\frac{\partial}{\partial x_j}\left(\frac{\partial v_i}{\partial x_j}\right) + \frac{1}{F_r^2}g_i, \qquad (9)$$

$$\frac{\partial v_i}{\partial x_i} = 0, \qquad (10)$$

para $i = 1, 2, 3$. Em (9-10) $v_i$ são as componentes do campo de velocidade do escoamento e $p$ é o respectivo campo de pressão. As quantidades $R_e = LU_o/v$ e $F = U_o/\sqrt{L|\mathbf{g}|}$ são, respectivamente, os números de Reynolds e Froude, enquanto $g_i$ são as componentes da força gravitacional.

### Salto hidráulico circular

Como uma primeira aplicação do esquema ADBQUICKEST considera-se o fenômeno do salto hidráulico circular (ELLEGAARD et al., 1998). Para descrever este tipo de escoamento utilizou-se o sistema de simulação Freeflow 3D (CASTELO FILHO et al., 2000), equipado com o esquema ADBQUICKEST. Foram considerados os seguintes dados: número de Reynolds $R_e = 150$, malha de $100 \times 100 \times 20$ células computacionais e tempo final de simulação $t = 0.5$s. A solução numérica está apresentada na Figura 5, onde pode-se constatar que o fenômeno físico foi simulado com sucesso.

**Figura 5 -** Simulação numérica tridimensional do salto hidráulico circular.

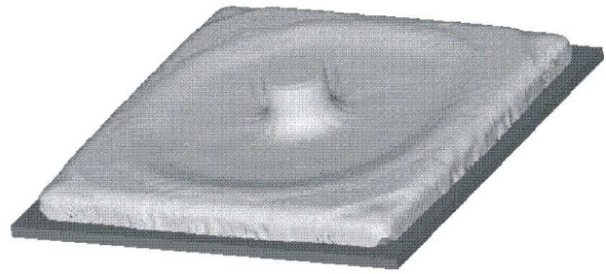

### Dobras em jatos livres

A segunda aplicação do esquema ADBQUICKEST simula o fenômeno de dobras (instabilidades físicas) em jatos livres a baixo número de Reynolds (TOMÉ et al., 2001). O propósito aqui é mostrar que o esquema ADBQUICKEST é apropriado também para problemas altamente viscosos. O código Freeflow 3D, incrementado com o esquema ADBQUICKEST, foi empregado na solução deste problema, com $R_e = 0.25$ numa malha de $100 \times 100 \times 100$ células computacionais e tempo final $t = 1.5$s. Como pode ser observado na Figura 6, o fenômeno físico de dobras foi simulado com sucesso nos casos planar e circular.

**Figura 6 -** Simulação numérica tridimensional de jatos oscilantes, casos planar e circular.

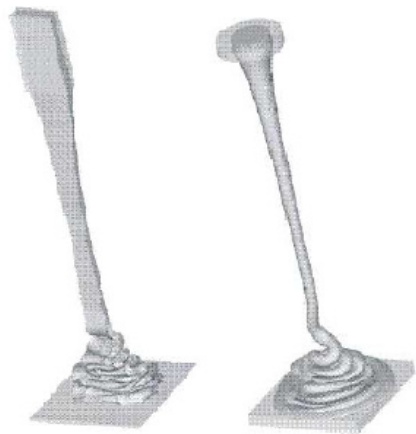

## Conclusões

O artigo forneceu uma avaliação computacional dos esquemas *upwind* populares de alta resolução WACEB, CUBISTA e ADBQUICKEST para resolver problemas não lineares em dinâmica dos fluidos. Em particular, foram simulados problemas de Riemann para as equações de Burgers, de Euler e de águas rasas, e, nestes problemas não lineares, o esquema ADBQUICKEST mostrou ser superior. Como aplicação, dois escoamentos de fluidos incompressíveis tridimensionais com superfícies livres móveis foram simulados com o esquema ADBQUICKEST, onde foi constatado que este constitui uma ferramenta efetiva para a simulação de fenômenos físicos altamente complexos, tais como o salto hidráulico circular e as instabilidades físicas em problemas a baixos números de Reynolds. Para o futuro, os autores planejam investigar o desempenho do esquema ADBQUICKEST em escoamentos de fluidos compressíveis sobre aerofólios e problemas envolvendo turbulência nos fluidos viscoelásticos.

## Agradecimentos

## Referências